\documentclass[10pt]{article}
\usepackage{amsmath,amsthm,amsfonts}
\usepackage{graphicx}
\usepackage{longtable}
\usepackage{float}
\usepackage{url}
\newdimen\epsfxsize
\newdimen\epsfysize
\def\qed{\vrule height5pt width3pt depth.5pt}

\theoremstyle{plain}
\newtheorem{thm}{Theorem}[section]
\newtheorem{cor}[thm]{Corollary}
\newtheorem{lem}[thm]{Lemma}

\newtheorem{conj}[thm]{Conjecture}
\newtheorem{exa}{Example}[section]

\newtheorem{defn}{Definition}[section]

\newtheorem{rem}{Remark}[section]

\begin{document}

\title{Parity Biquandles}


\author{Aaron Kaestner \\
University of Illinois at Chicago\\
akaestne@math.uic.edu \\
\\
Louis H. Kauffman \\
University of Illinois at Chicago \\
kauffman@uic.edu }

\maketitle

\begin{abstract} We use crossing parity to construct a generalization of biquandles for virtual knots which we call Parity Biquandles.  These structures include all biquandles as a standard example referred to as the even parity biquandle.  Additionally, we find all Parity Biquandles arising from the Alexander Biquandle and Quaternionic Biquandles. For a particular construction named the z-Parity Alexander Biquandle we show that the associated polynomial yields a lower bound on the number of odd crossings as well as the total number of real crossings and virtual crossings for the virtual knot. Moreover we extend this construction to links to obtain a lower bound on the number of crossings between components of a virtual link.
\end{abstract}

\section{Introduction}

\subsection{Virtual Knots and Biquandles}

In \cite{VKT} Kauffman introduced virtual knots and links as a natural extension of classical knot theory.  Virtual Knot Theory can be though of both as 1), equivalent classes of an embedded closed curve in a thickened surface $S_g \times I $ (possibly non-orientable) up to isotopy and handle stabilization on the surface and 2) the completion of the oriented Gauss codes (i.e. an arbitrary Gauss code corresponds to a virtual knot while not every Gauss code corresponds to a classical knot.)\\

Invariants for virtual knots arising from the analysis of chord diagrams were introduced in \cite{VKT} and further explored by Goussarov, Polyak and Viro in \cite{GPV}. Biquandles have a rich history in virtual knot theory including work by Sawollek (\cite{Sawollek1}), Nelson (\cite{NelsonBiracks}), Fenn, Kauffman and Jordan-Santana(\cite{FJK}) Kauffman and Manturov (\cite{BQforVirts}), Kauffman and Hrencecin (\cite{VirtBQs}), Kauffman and Radford (\cite{KR1}), and Bartholomew and Fenn (\cite{BFSmall}, \cite{Quaternionic}). Similarly, virtual knot invariants arising from an analysis of parity have previously been constructed by Kauffman (\cite{SelfLinking}, \cite{IntroVKT}), Manturov (\cite{ManturovParity}), Turaev (\cite{Turaev}) and Dye (\cite{DyeLinkingNumber}).  Our approach to parity was inspired by Manturov's philosophy of parity (\cite{ManturovIlyutkoNikonov},\cite{ManturovKrylov},\cite{ManturovParity}) and the construction of the parity bracket polynomial (\cite{ManturovPKT}, \cite{IntroVKT}).\\

\begin{figure}[H]
\centering
    \includegraphics[width=0.5\textwidth]{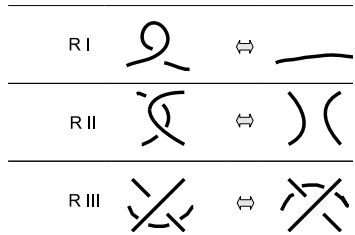}
\caption{Reidemeister Moves}
\label{fig:RMs}
\end{figure}

\begin{figure}[H]
\centering
    \includegraphics[width=0.5\textwidth]{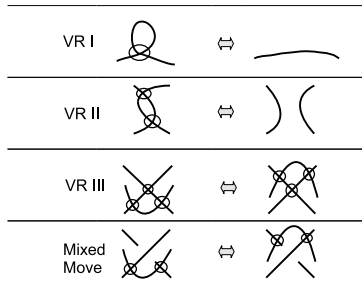}
\caption{Virtual Reidemeister Moves}
\label{fig:VRMs}
\end{figure}

We recall in Figure \ref{fig:RMs} the Reidemeister Moves and their corresponding flat moves on chord diagrams. Figure \ref{fig:VRMs} displays the additional Virtual Reidemeister Moves, note these have no affect on the chord diagram.\\

Following \cite{NelsonBiracks} and \cite{KR1} we recall the definition of a Biquandle.\\

\begin{defn} A biquandle $(X,B)$ is a set $X$ and a map $B : X\times X \rightarrow X\times X$ which satisfies the following conditions:
\begin{enumerate}
	\item $B$ is invertible, i.e there exists a map $B^{-1}: X\times X \rightarrow X\times X$ satisfying $B\circ B^{-1} = Id_{X\times X} = B^{-1}\circ B$,\\
	\item For all $a,b \in X$ there exists $x \in X$ such that \\ $x = B_{2}^{-1}(a,B_{2}(b,x))$, $a = B_{1}(b,x)$ and $b = B_{1}^{-1}(a, B_{2}(b,x))$ \\
	 For all $a,b \in X$ there exists $x \in X$ such that \\ $x = B_{1}(B_{1}^{-1}(x,b),a)$, $a = B_{2}^{-1}(x,b)$ and $b = B_{2}(B_{1}^{-1}(x,b),a)$ \\
	\item $B$ satisfies the set-theoretic Yang-Baxter equation $(B\times Id) \circ (Id\times B) \circ (B\times Id) = (Id\times B) \circ (B\times Id) \circ (Id\times B)$\\
	\item Given $a \in X$ there exists $x \in X$ such that $a = B_{1}(a,x)$ and $x = B_{2}(a,x)$  \\ Given $a \in X$ there exists $x \in X$ such that $a = B_{1}^{-1}(a,x)$ and $x = B_{2}^{-1}(a,x)$  \\
\end{enumerate}
\end{defn}

\begin{figure}[H]
\centering
    \includegraphics[height=3cm]{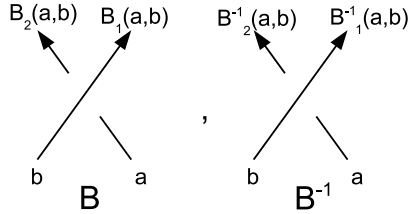}
\caption{Diagrammatic Representation of the Biquandle}
\label{fig:BqDiagramatic}
\end{figure}

Diagrammatically $B$ and $B^{-1}$ corresponds to a crossing as in Figure \ref{fig:BqDiagramatic}. Reinterpreting the above definition in this diagrammatic form we see that the Axioms 1 and 2 for $B$ are equivalent to the same-oriented and opposite-oriented Reidemeister II Moves, Axiom 4 corresponds to a Reidemeister I Move and Axiom 3 corresponds to a same-oriented, positive crossing Reidemeister III Move. It is a simple exercise (\cite{KnotsandPhysics}) to show that this is enough to ensure invariance under all remaining oriented Reidemeister Moves.  \\

Given a knot $K$, the biquandle of the knot $K$, $BQ(K)$, is the non-associative algebra generated by the arcs in any planar diagrams of $K$ and relations given by the map $B$.  \\

\begin{lem} $BQ(K)$ is an invariant of the virtual knot $K$.
\end{lem}

\begin{rem} Those familiar with the subject will note that the removal of Axiom 4 from the above list gives the definition of a \textit{birack}. This omission, along with the following section, yields the appropriate definition of \textit{parity birack}.  We will not discuss parity biracks further other than to remark that, just as every biquandle is a birack, every parity biquandle is a parity birack.\\
\end{rem}

Some common examples of biquandles are the Generalized Alexander Biquandle \cite{KR1}, \cite{Sawollek1} and the Quaternionic Biquandles with integral coefficients \cite{Quaternionic}.\\

\begin{figure}[H]
\centering
    \includegraphics[height=3cm]{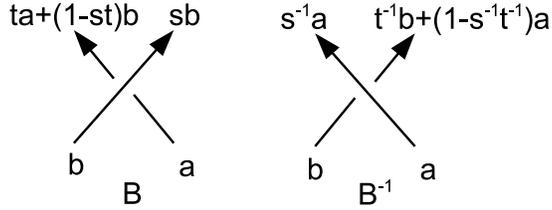}
\caption{Generalized Alexander Biquandle}
\label{fig:GAB}
\end{figure}

The Generalized Alexander Biquandle is defined by the diagram in Figure \ref{fig:GAB} where $a,b \in X$, where $s$ and $t$ are commuting variables in the ground ring, and results in a  $\mathbb{Z}\left[ s^{\pm 1},t^{\pm 1} \right]$-module. The following example shows how to use this definition to arrive at the Sawollek Polynomial (\cite{KR1}, \cite{Sawollek1}) , a Laurent Polynomial in $\mathbb{Z}\left[ s^{\pm 1},t^{\pm 1} \right]$. Note this polynomial is unique up to a multiple of $t^{\pm 1}$.\\

\begin{exa}

Consider the 3-Crossing Knot $3.1$  (our naming conventions follow Jeremy Green's Knot Tables \cite{JerGreen}) in Figure \ref{fig:3.1}.

\begin{figure}[H]
\centering
    \includegraphics[height=3cm]{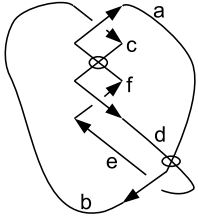}
\caption{Virtual Knot 3.1}
\label{fig:3.1}
\end{figure}

Following the convention of Figure \ref{fig:GAB} we obtain the following system of equations:

\begin{center}
	\[
\begin{array}{rcl}
	a & = & s^{-1}f\\
	b & = & s^{-1}a\\
	c & = & t^{-1}b + (1-s^{-1}t^{-1})f\\
	d & = & sc\\
	e & = & t^{-1}d + (1-s^{-1}t^{-1})a\\
	f & = & te + (1-st)
\end{array}
\]
\end{center}

Or equivalently:\\

\begin{center}
\[
	\begin{array}{lcl}
	-a + s^{-1}f & = &  0\\
	s^{-1}a - b & = &  0\\
	t^{-1}b - c + (1-s^{-1}t^{-1})f & = &  0\\
	sc - d & = &  0\\
	(1-s^{-1}t^{-1})a + t^{-1}d - e & = &  0 \\
	(1-st)c +te -f & = &  0
\end{array}
\]
\end{center}

Fixing the basis $\left\{ a, b, c, d, e, f\right\}$ of $X^{\times 6}$ we obtain the matrix:\\
\[	\begin{pmatrix}
	-1 & 0 & 0 &  0 & 0 & s^{-1}\\
	s^{-1} & -1 & 0 &  0 & 0 & 0 \\
	0 & t^{-1} & -1 & 0 & 0 & (1-s^{-1}t^{-1})\\
	 0 & 0 & s & -1 & 0 & 0  \\
	(1-s^{-1}t^{-1}) & 0 & 0 & t^{-1} & -1 & 0 \\
	0 & 0 & (1-st) &  0 & t &-1 \\
	\end{pmatrix}
\]

Taking the determinant and multiplying by $(-1)^{wr(K)}$, where $wr(K)$ = writhe($K$) = (\# positive crossings) - (\# negative crossings), we find, up to multiples of $s^{n}t^{m}$, $n, m \in \mathbb{Z}$, the Sawollek Polynomial of virtual knot $3.1$ is
	\[ \frac{1-\frac{1}{s^2}}{t}+\frac{1}{s^2}+\left( s-\frac{1}{s} \right) t-s+\frac{1}{s}-1
\]
\end{exa}

For a more systematic description of the matrix construction see \cite{Sawollek1}.  It should be noted that the Sawollek polynomial and the generalizations presented later in this paper are well-defined following the proof given in \cite{Quaternionic} and in the spirit of \cite{CrowellFox}.  When working over a gcd-ring, including a polynomial ring over $\mathbb{Z}$,  the determinant of the presentation matrix generates a principle ideal and is an invariant of the knot \cite{Quaternionic}.  Recall that for a classical knot one of the relations in the matrix above will always be a consequence of the others, hence the Sawollek polynomial will be identically zero on classical knots.  \\

\begin{figure}[H]
\centering
    \includegraphics[height=3cm]{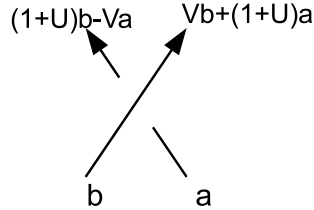}
\caption{Quaternionic Biquandle with Integral Coefficients}
\label{fig:QB}
\end{figure}

As described in \cite{Quaternionic} the Quaternionic Biquandles with integral coefficients are a defined as in Figure \ref{fig:QB} where $U, V  \in \left\{ \pm i, \pm j, \pm k \right\}, U \bot V$.\\

\subsection{Parity and Virtual Knots}

Given a diagram $D$ for a knot $K$ label each crossing uniquely 1 through $n$, where $n$ is the total number of crossings in $D$.  Let $P$ an arbitrary base-point on the knot. Starting at $P$ and following the orientation of the knot we can construct a sequence of length $2n$ with terms corresponding to each crossing we encounter. Each term is a 3-tuple of the form ($O$/$U$, Crossing Number, $\pm$) where $O$ or $U$ corresponds to an over or under-crossing respectively and $\pm$ corresponds to the sign of the crossing.  The resulting code is referred to as the (signed, oriented) {\itshape Gauss Code} for the diagram $D$ of the knot $K$. \\

\begin{figure}[H]
\centering
    \includegraphics[height=3cm]{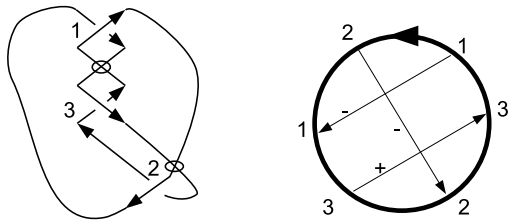}
\caption{Chord Diagram for Virtual Knot 3.1}
\label{fig:3_1CD}
\end{figure}

The Gauss code can be represented diagrammatically as follows.  Given a circle (often referred to as the {\itshape core circle}) place upon it in a counterclockwise fashion $2n$ points where each point is labeled by a crossing name (an integer between 1 and n) in the cyclic order corresponding to the Gauss code.  Between the two occurrences of a crossing on the core circle, place an signed, oriented chord where the sign corresponds to the crossing sign and the orientation goes from the over crossing to the under crossing. We call this the {\itshape Chord Diagram} for $D$. (\cite{GPV}, \cite{VKT}) For example, the knot $3.1$ in Figure \ref{fig:3.1} has Gauss Code ``$01-,02-,U1-,O3+,U2-,U3+$'' and chord diagram as in Figure \ref{fig:3_1CD}.  \\

\begin{defn}
Given a (virtual) knot \textbf{K} we can label each crossing as even or odd in the following manner. For each crossing $v$ locate the 2 occurrences of $v$ in the Gauss code for \textbf{K}. If the number of crossing labels between the two occurrences of $v$ is even then label the crossing even.  Else it is labeled odd.\\
\end{defn}

\begin{rem} This parity is well-defined for a 1-component links (i.e. knots) as the number of crossing labels in the Gauss code is $2n$ where $n$ is the number of crossings.\\
\end{rem}

It is important to notice how parity behaves under the classical Reidemeister moves, recalling that virtual Reidemeister moves do not change the Gauss code or chord diagram and thus do not affect parity. \\

\begin{figure}[H]
\centering
    \includegraphics[width=0.5\textwidth]{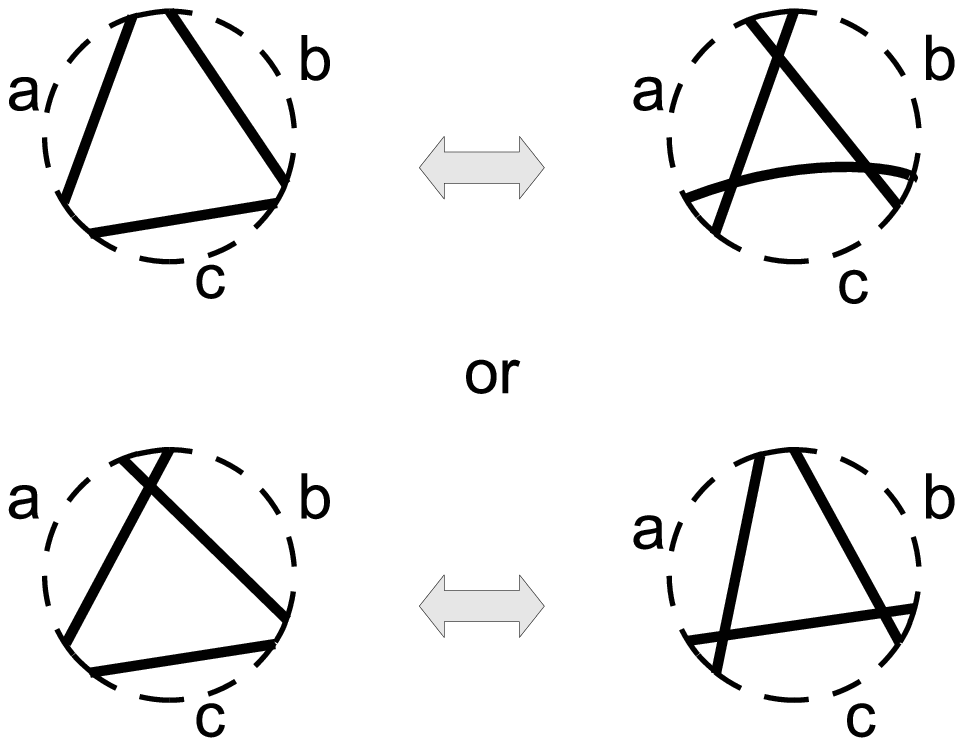}
\caption{}
\label{fig:CDR3}
\end{figure}

\begin{itemize}
	\item \textbf{Reidemeister I}\\
	A first Reidemeister move is always even, as is shown in Figure \ref{fig:CDR1}
	
	\begin{figure}[H]
\centering
    \includegraphics[width=0.5\textwidth]{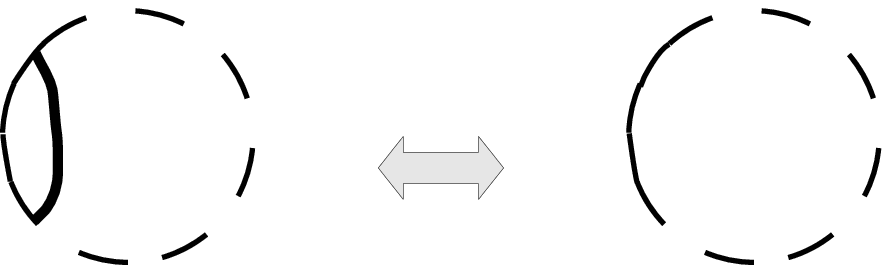}
\caption{}
\label{fig:CDR1}
\end{figure}
	
	\item \textbf{Reidemeister II}\\
	The two crossings involved in a second Reidemeister move are either both even or both odd.  To see this, note that in Figure \ref{fig:CDR2}  if the number of crossings before the second Reidemeister move is $n+2$ and $a$ and $b$ denote the number of markings on the core circle as labeled in the figure then $a + b = 2n$ is even. Hence either $a$ and $b$ are both even or both odd.
	
	\begin{figure}[H]
\centering
    \includegraphics[width=0.5\textwidth]{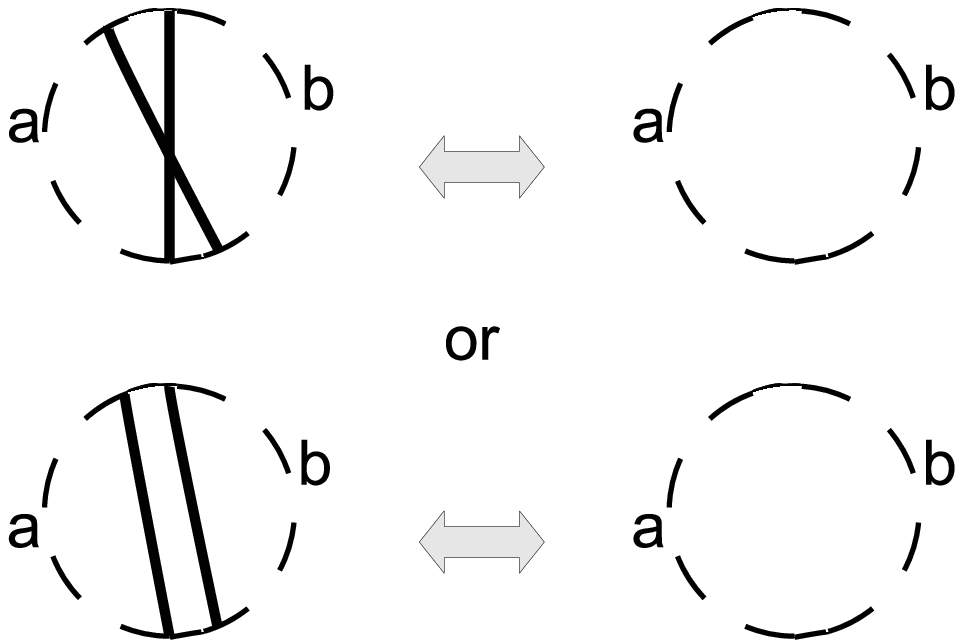}
\caption{}
\label{fig:CDR2}
\end{figure}

\item \textbf{Reidemeister III}\\
	In a third Reidemeister move either all crossings are even or two are even and one is odd. To see this note that in Figure \ref{fig:CDR3}  if the number of crossings not involved in the third Reidemeister move is $n$ and $a, b$ and $c$ denote the number of markings on the core circle as labeled in the figure then $a + b + c= 2n$ is even. Hence either $a, b$ and $c$ are all even or two are even and one is odd.
\end{itemize}

Thus we can generalize the definition of the biquandle by constructing separate maps for the odd and even crossings.

\section{Parity Biquandles}

\begin{defn} A Parity Biquandle $(X,B,P)$ is biquandle $(X,B)$ and a map $P: X\times X \rightarrow X\times X$ which satisfies

\begin{enumerate}
	\item $P$ is invertible, i.e there exists a map $P^{-1}: X\times X \rightarrow X\times X$ satisfying $P\circ P^{-1} = Id_{X\times X} = P^{-1}\circ P$,\\
	\item
	For all $a,b \in X$ there exists $x \in X$ such that \\ $x = P_{2}^{-1}(a,P_{2}(b,x))$, $a = P_{1}(b,x)$ and $b = P_{1}^{-1}(a,P_{2}(b,x))$ \\
	For all $a,b \in X$ there exists $x \in X$ such that \\ $x = P_{1}(P_{1}^{-1}(x,b),a)$, $a = P_{2}^{-1}(x,b)$ and $b = P_{2}(P_{1}^{-1}(x,b),a)$ \\
	\item $B$ and $P$ satisfy the set-theoretic Yang-Baxter equations \\
	$(P\times Id) \circ (Id\times P) \circ (B\times Id) = (Id\times B) \circ (P\times Id) \circ (Id\times P)$ \\
	$(P\times Id) \circ (Id\times B) \circ (P\times Id) = (Id\times P) \circ (B\times Id) \circ (Id\times P)$ \\
	$(B\times Id) \circ (Id\times P) \circ (P\times Id) = (Id\times P) \circ (P\times Id) \circ (Id\times B)$ \\
\end{enumerate}
\end{defn}

\begin{defn} Given a biquandle $(X,B)$ the even parity biquandle of $(X,B)$ is the parity biquandle $(X,B,B)$.
\end{defn}

Given a knot $K$ and diagram $D$, the parity biquandle of the knot $K$, $PBQ(K)$, is the non-associative algebra generated by the arcs in $D$ and relations given by applying the maps $B$ at even crossings of $D$ and $P$ at odd crossings of $D$.  \\

\begin{lem} $PBQ(K)$ is an invariant of the virtual knot $K$.
\end{lem}

\subsection{The Parity Alexander Biquandle}

Given that $B : X\times X \rightarrow X\times X$ as described in Figure \ref{fig:GAB} is linear we can represent $B$ by the matrix
	\[
	\left[ \begin{array}{cc}
	0 & s \\ t & 1 - st
	\end{array} \right] \]

Representing $P$ by a $2 \time 2$ matrix, we utilized the linear algebra functionality of Mathematica to determine the following possible values for $P$.

\begin{enumerate}
	\item  Even Parity Alexander Biquandle\[ P_1 = B =
		\left[ \begin{array}{cc}
	0 & s \\ t & 1 - st
	\end{array} \right] \]

	\item \[ P_2 =
		\left[ \begin{array}{cc}
	0 & s \\ t &  st-1
 \end{array} \right] \]	

	\item z-Parity Alexander Biquandle \[ P_3 =
  	\left[ \begin{array}{cc}
	0 & z \\ z^{-1} &  0
\end{array} \right] \]  \\
\end{enumerate}

or diagrammatically:\\

\begin{figure}[H]
\centering
    \includegraphics[width=0.5\textwidth]{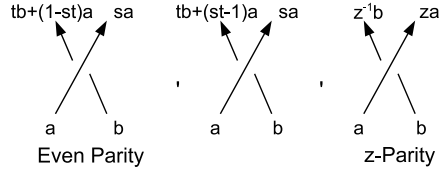}
\caption{Diagrammatic Representations for P in the Parity Alexander Biquandle}
\label{fig:PABqs}
\end{figure}

Thus $(X,B,P_1)$, $(X,B,P_2)$ and $(X,B,P_3)$ are each parity biquandles. Note $(X,B,P_1)$ and $(X,B,P_2)$ generate  $\mathbb{Z}\left[ s^{\pm 1},t^{\pm 1} \right]$-modules  while $(X,B,P_3)$, the z-Parity Alexander Biquandle, generates a  $\mathbb{Z}\left[ s^{\pm 1},t^{\pm 1},z^{\pm 1}  \right]$-module.\\

Note that $P_2  = \begin{bmatrix} 0 & s \\ t & 1 - st \end{bmatrix} \times \begin{bmatrix} 1 & 2(s-t^{-1}) \\ 0 & 1 \end{bmatrix}$.
Although the polynomial invariant induced by $P_2$ appears distinct from the Sawollek polynomial, we have yet to find any computational benefit resulting from its calculation. $P_3$ is a different matter.  Namely we have the following theorems:\\

\begin{thm}
\begin{enumerate}
	\item If the polynomial associated to any Parity Alexander Biquandle for a virtual knot $K$ is nonzero, then $K$ is nonclassical.
	\item If the z-Parity Alexander Polynomial for a virtual knot $K$ is unequal to the Sawollek Polynomial for $K$ then any diagram of $K$ contains an odd crossing.\\
\end{enumerate}
\end{thm}

\textbf{Proof:}
\begin{enumerate}
	\item Suppose $K$ is (equivalent to) a classical knot. By the Jordan Curve Theorem $K$ is equivalent to a knot with no odd crossings.  Thus any Parity Alexander Biquandle for $K$ is equivalent to the Generalized Alexander Biquandle for $K$. Thus the respective polynomial is equivalent to the Sawollek polynomial which is identically $0$ on classical knots. (See Theorem 3 in \cite{Sawollek1} or \cite{JaegKaufSal})\\
	\item Similarly, if $K$ has a diagram $D$ with no odd crossings, then the z-Parity Alexander Biquandle of $D$ is equivalent to the Generalized Alexander Biquandle of $D$. Since the Parity Biquandle is invariant under the Reidemeister moves we have the z-Parity Alexander Polynomial of $K$ is equal to the Sawollek Polynomial for $K$.\\
\end{enumerate}
\qed

Moreover, the z-Parity Alexander Polynomial provides a lower bound on the minimum number of odd crossings in a virtual knot.\\

\begin{thm}\label{thm:main}
Given a virtual knot $K$, let $n_{o}$ be the minimum number of odd crossings in any diagram of $K$, and suppose $z^{e_{max}}$ and $z^{e_{min}}$ are, respectively, the highest and lowest powers of $z$ appearing in the z-Parity Alexander Polynomial of $K$, and set $e = \max{(\left| e_{max} \right|,\left| e_{min} \right|)}$ then

$\begin{cases} e \leq n_{o}, & \mbox{if }  e \mbox{ is even} \\ (e + 1) \leq n_{o}, & \mbox{if } e \mbox{ is odd}  \end{cases}$\\

\end{thm}

\textbf{Proof:}
Suppose $D$ is a diagram for $K$ with a minimal number of odd crossings and let $n$ be the number of odd crossings in $D$.  Then the matrix of relations contains $n$ entries of value $z$ (and $z^{-1}$). Thus the highest and lowest power of $z$ in any term of the determinant is $\pm n$. This gives the inequality $e \leq n_{o}$. Since the number of odd crossings in any knot is always even (Prop. 1.2 \cite{DyeLinkingNumber}) we get the theorem.
\begin{center}
\qed \\
\end{center}

\begin{cor}\label{cor:maincor}
Given a virtual knot $K$ let $n$ be the minimum number of real (non-virtual) crossings in any diagram of $K$ and define $e$ as in the previous theorem. If $e > 0$ then $ (e + 1) \leq n$.
\end{cor}

\textbf{Proof:}
Let $D$ be any diagram for $K$. Suppose for contradiction $e > 0$ and  $D$  has no even crossings.  Then every relation is of the form $a = z ^{\pm 1} b$ where $a$ and $b$ are consecutive arc labels of $D$. Moreover, starting at any arc and traversing $D$ creates a cycle of relations of the above form with $n$ occurrences of $z$ and $z^{-1}$. Hence the z-Parity Alexander Biquandle of $K$ is trivial and thus $e = 0$.
\begin{center}
\qed \\
\end{center}

\begin{exa} Following the same procedure as earlier with knot $3.1$ we see it has z-Parity Alexander Polynomial \[ \frac{\frac{1}{s t}-1}{z^2}-\frac{1}{s t}+1 \] Hence the diagram in Figure \ref{fig:3.1} is minimal for virtual knot $3.1$ in the sense that it contains the minimum number of odd crossings and the minimum number of total crossing for any diagram of the knot.\\
\end{exa}

\begin{figure}[H]
\centering
    \includegraphics[height=3cm]{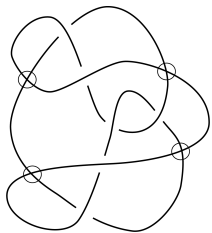}
\caption{Knot 6.32008}
\label{fig:6_32008}
\end{figure}

\begin{figure}[H]
\centering
    \includegraphics[height=3cm]{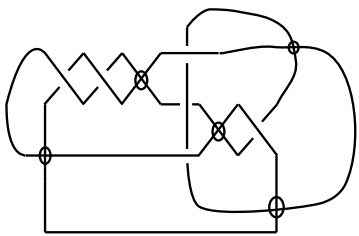}
\caption{Knot 6.73583}
\label{fig:6_73583}
\end{figure}

Using Jeremy Green's tables \cite{JerGreen} we have calculated the Sawollek Polynomial and the z-Parity Alexander Polynomial for knots with at most 6 real crossings.  The knots in Figure \ref{fig:6_32008} and Figure \ref{fig:6_73583} are special in that they are not distinguished from the unknot via the Sawollek Polynomial, z-Parity Sawollek Polynomial, Arrow Polynomial and Parity Arrow Polynomial \cite{ParityKhoHoArrowCat}. Knot $6.32008$ has 4 odd crossings while Knot $6.73583$ has no odd crossings and both knots are trivial as flats. Using a 2-cable Jones Polynomial calculator adapted from Dror Bar-Natan's "faster" Jones Polynomial Calculator \cite{DrorJones} we have been able to distinguish each of these knots from one-another and from the unknot.\\

See the Appendix for calculations on knots with at most 4 real crossings following the conventions of \cite{JerGreen} along with their Sawollek Polynomials. Note that of the 19 knots with Sawollek polynomial equal 0, 3 are detected (nonzero) by z-Parity. Similarly, of the 54 knots with z-Parity polynomial equal 0, 38 are detected by Sawollek.\\

While investigating computations for the z-Parity Alexander polynomial we have verified the following conjecture on virtual knots with less than 6 real crossings.\\

\begin{conj}
Given a virtual knot $K$ let $n$ be the minimum number of real (non-virtual) crossings in any diagram of $K$ and suppose $z^{e_{max}}$ and $z^{e_{min}}$ are, respectively, the highest and lowest powers of $z$ appearing in the z-Parity Alexander Polynomial of $K$. Then $(e_{max} - e_{min}) \leq n $\\
\end{conj}

The lower bound in this conjecture rarely appears to be tight.  The following example is one of five knots with 4-crossings where the bound equals the minimum real crossing number.\\

\begin{exa}
Knot 4.96, given by Gauss Code
	\[`` O1-,O2-,U3+,U1-,O4-,U2-,O3+,U4- ''\]
	and having 2 odd crossings, has z-Parity Alexander Polynomial
	\[z^2 \left(\frac{1}{s t}-\frac{1}{s^2 t^2}\right)+\frac{\frac{1}{s t}-\frac{1}{s^2 t^2}}{z^2}+\frac{2}{s^2 t^2}-\frac{2}{s t}\]
\end{exa}

\subsection{Parity Quaternionic Biquandles}

In the same fashion as the Parity Alexander Biquandle we utilized Mathematica along with the matrix representation to determine the following values for $P$, when $(X,B)$ is a Quaternionic Biquandle with integral coefficients.\\

\begin{figure}[H]
\centering
    \includegraphics[height=6cm]{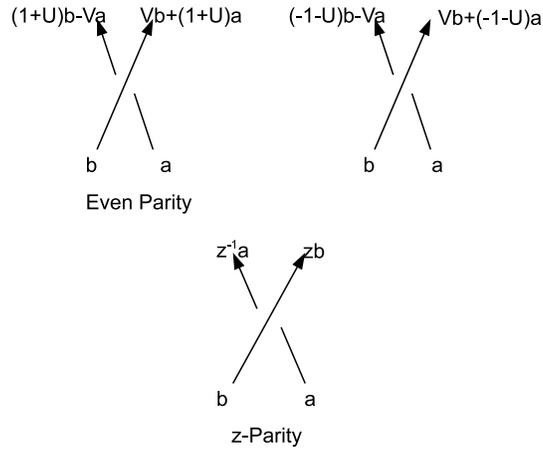}
\caption{Diagrammatic values for P for the Parity Quaternionic Biquandle}
\label{fig:PQBqs}
\end{figure}

\begin{figure}[H]
\centering
    \includegraphics[height=2.5cm]{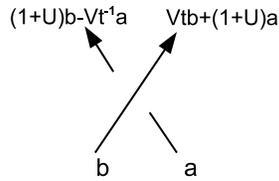}
\caption{Diagrammatic Representation for B in the Quaternionic Biquandle with integral coefficients after a change of basis}
\label{fig:BasisChange}
\end{figure}

To create an polynomial invariant from the quaternionic biquandle we follow the construction in \cite{Quaternionic}. We first perform a change of basis on the map $B$ which corresponds to extending the ground ring by commuting variables $t, t^{-1}$ .  This can be represented diagrammatically as in Figure \ref{fig:BasisChange}.  The construction follows analogously to the Sawollek polynomial. However, before taking the determinant we replace each element of the presentation matrix with its corresponding $SU \left( 2 \right)$ matrix representation. For an n-crossing knot this produces a $4n \times 4n$ matrix over $\mathbb{C}$ whose determinant, called the Study Determinant in \cite{Quaternionic}, is an invariant of the knot.\\

 \begin{exa} Setting $U = i$ and $V = j$ in Figure \ref{fig:QB} the virtual knot $3.1$ in Figure \ref{fig:3.1} has z-Parity Quaternionic polynomial\\
	\[2 z^4+\frac{2}{z^4}-4 z^2-\frac{4}{z^2}+4\]
 \end{exa}

\subsection{Link Parity Biquandles}

\begin{figure}[H]
\centering
    \includegraphics[width=0.5\textwidth]{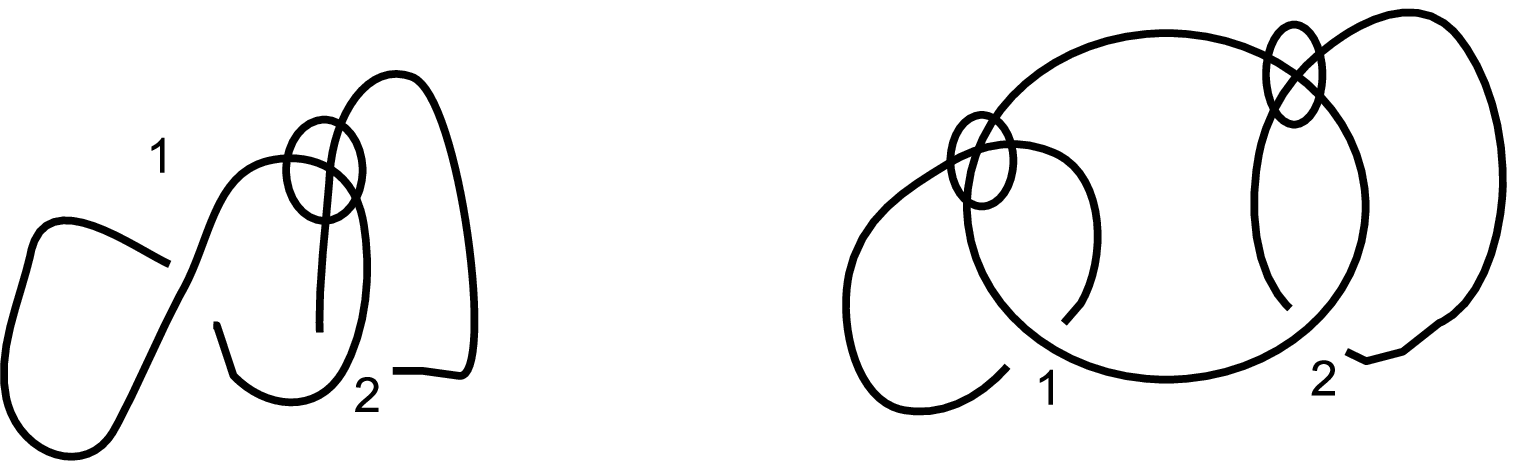}
\caption{}
\label{fig:paritylinks1}
\end{figure}

One should note that our definition of even and odd parity does not naturally extend to links (2 or more components). For example, the links in Figure \ref{fig:paritylinks1} illustrate some of the difficulty in the natural extension.\\

Omitting signs, the left link in \ref{fig:paritylinks1} has Gauss code ``$O1,U1,O2; U2$'' while the other has Gauss code ``$U1; O1,O2; U2$''. In the first of these Crossing 1 is both even and odd in the first component while Crossing 1 is either even or odd depending upon whether you examine the first or second link component.\\

We may circumvent this pitfall by defining even and odd for self-crossings based on the parity of self-crossing in each component while labeling crossings shared by 2 components as link crossings.

\begin{figure}[H]
\centering
    \includegraphics[height=3cm]{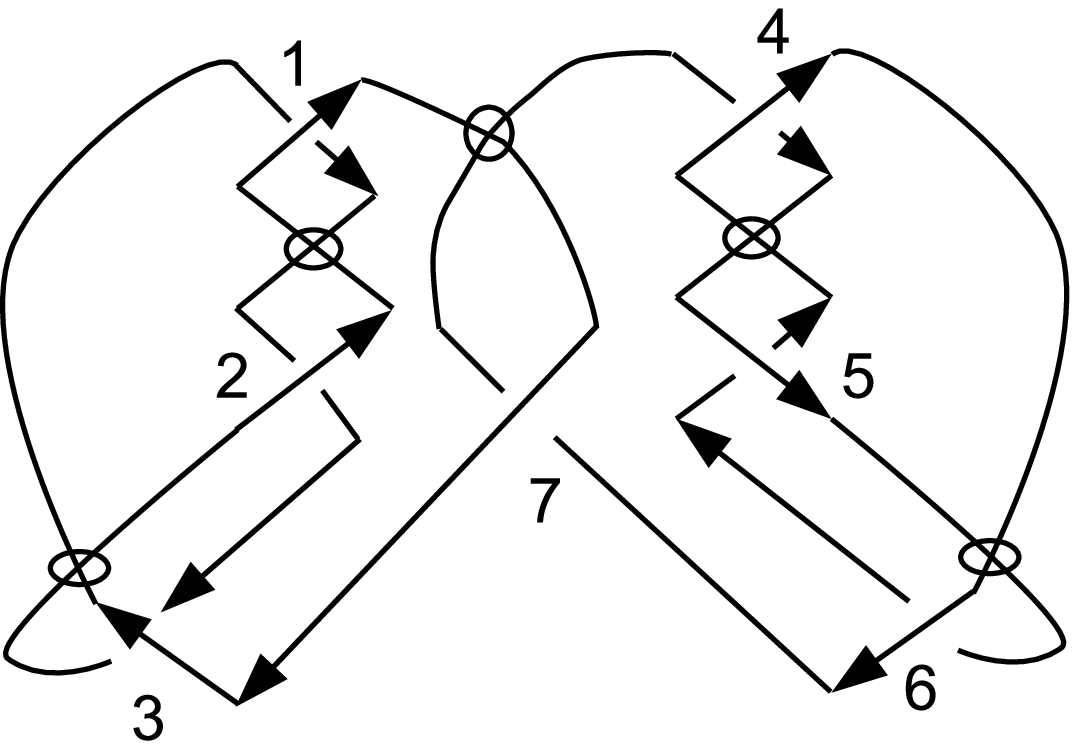}
\caption{}
\label{fig:LinkParityExample}
\end{figure}

\begin{exa}The link in Figure \ref{fig:LinkParityExample} has Gauss Code \[``01,O7,03,U1,U2,U3,O2; U4,O5,U6,U5,O4,O6,U7'' \]  Crossings 1, 2, 4 and 5 are odd, crossings 3, and 6 are even and crossing 7 is a link crossing.
\end{exa}

As we did with odd crossings, we investigate the invariance of crossings between links to provide the framework for generalizing the Parity Biquandle to the Link Parity Biquandle. We will call a crossing where both arcs involved are in one link component an self-crossing while a crossing whose arcs are in separate components a link crossing.\\

\begin{itemize}
	\item \textbf{Reidemeister I}\\
	In a Reidemeister I move only a single link component is involved, as is shown in Figure \ref{fig:RMs}.\\
	
	\item \textbf{Reidemeister II}\\
	The two strands involved in a second Reidemeister move are either both in the same component or each in a  different component as is shown in Figure \ref{fig:RMs}.  Thus either both crossings above are self-crossings or both crossings are link crossings.\\

\item \textbf{Reidemeister III}\\
	In a third Reidemeister move either all strands involved are in one component, or two in one component and one in another or all three in separate components in Figure \ref{fig:RMs}. Thus either all crossings are self-crossings, or there is one self-crossing and two link crossings or three link crossings respectively.\\
\end{itemize}

\begin{defn} A Link Parity Biquandle $(X,B,P,L)$ is biquandle $(X,B,P)$ and a map $L: X\times X \rightarrow X\times X$ which satisfies

\begin{enumerate}
	\item $L$ is invertible, i.e there exists a map $L^{-1}: X\times X \rightarrow X\times X$ satisfying $P\circ P^{-1} = Id_{X\times X} = P^{-1}\circ P$,\\
	\item
	For all $a,b \in X$ there exists $x \in X$ such that \\ $x = L_{2}^{-1}(a,L_{2}(b,x))$, $a = L_{1}(b,x)$ and $b = L_{1}^{-1}(a,L_{2}(b,x))$ \\
	For all $a,b \in X$ there exists $x \in X$ such that \\ $x = L_{1}(L_{1}^{-1}(x,b),a)$, $a = L_{2}^{-1}(x,b)$ and $b = L_{2}(L_{1}^{-1}(x,b),a)$ \\
	\item $B$, $P$,and $L$ satisfy the set-theoretic Yang-Baxter equations \\
	$(L\times Id) \circ (Id\times L) \circ (B\times Id) = (Id\times B) \circ (L\times Id) \circ (Id\times L)$ \\
$(L\times Id) \circ (Id\times B) \circ (L\times Id) = (Id\times L) \circ (B\times Id) \circ (Id\times L)$ \\	 $(B\times Id) \circ (Id\times L) \circ (L\times Id) = (Id\times L) \circ (L\times Id) \circ (Id\times B)$ \\	 $(L\times Id) \circ (Id\times L) \circ (P\times Id) = (Id\times P) \circ (L\times Id) \circ (Id\times L)$ \\	 $(L\times Id) \circ (Id\times P) \circ (L\times Id) = (Id\times L) \circ (P\times Id) \circ (Id\times L)$ \\	 $(P\times Id) \circ (Id\times L) \circ (L\times Id) = (Id\times L) \circ (L\times Id) \circ (Id\times P)$ \\
$(L\times Id) \circ (Id\times L) \circ (L\times Id) = (Id\times L) \circ (L\times Id) \circ (Id\times L)$ \\
\end{enumerate}
\end{defn}

Given a link $K$ and diagram $D$, the link parity biquandle of the knot $K$, $LPBQ(K)$, is the non-associative algebra generated by the arcs in $D$ and relations given by applying the maps $B$ at even crossings of $D$, $P$ at odd crossings of $D$ and $L$ at link crossings of $D$.  \\

\begin{lem} $LPBQ(K)$ is an invariant of the virtual link $K$.\\
\end{lem}

Furthermore, since at most one even or odd crossing can be involved in any Reidemeister move, we get the following generalization:\\

\begin{defn} A Generalized Link Parity Biquandle
\[(X,\{B_{\lambda}\}_{\lambda \in \Lambda},\{P_{\lambda}\}_{\lambda \in \Lambda},\{L_{\{\lambda, \rho\}}\})\]
is a family where for every $\lambda, \rho, \gamma \in \Lambda = \{1,\ldots,n \},  \lambda \neq \rho \neq \gamma$,  $(X, B_{\lambda}, P_{\lambda})$ is a parity biquandle, $(X, B_{\lambda}, P_{\lambda}, L_{\{\lambda, \rho\}})$ is a link parity biquandle and the maps satisfy the following condition:
\[(L_{\{\lambda, \rho\}} \times Id) \circ (Id\times L_{\{\lambda, \gamma\}}) \circ (L_{\{\rho, \gamma\}} \times Id) = (Id \times L_{\{\rho, \gamma\}}) \circ (L_{\{\lambda, \gamma \}}\times Id) \circ (Id\times L_{\{\lambda, \rho\}}) \]
\end{defn}

Given an n-component link $K$ with diagram $D$, and components labeled $1, \ldots, n$ the generalized link parity biquandle of the link $K$, $GPBQ(K)$, is the non-associative algebra generated by the arcs in $D$ and relations given by applying the map $B_{\lambda}$ at even crossings of component $\lambda$, the map $P_{\lambda}$ at odd crossings of component $\lambda$, and the map $L_{\{\lambda,\rho\}}$ at crossings between components $\lambda$ and $\rho$ for $\lambda \neq \rho \in \{1, \ldots, n\}$\\

\begin{lem} $GPBQ(K)$ is an invariant of the virtual link $K$.\\
\end{lem}

\subsection{The Generalized Link Parity Alexander Biquandle}

\begin{figure}[H]
\centering
    \includegraphics[height=2.5cm]{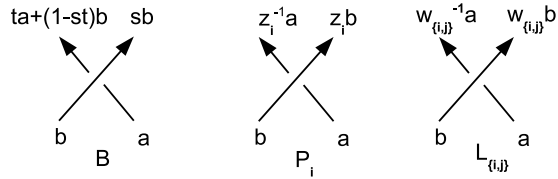}
\caption{Diagrammatic Relations for the Generalized Link Parity Alexander Biquandle}
\label{fig:GenParityAlexBQ}
\end{figure}

Suppose $(X,B)$ is the Generalized Alexander Biquandle described in Figure \ref{fig:GAB}.  We have shown that for a single component we may generalize via parity to the z-Parity Alexander Biquandle by applying the relation in Figure \ref{fig:PABqs} at odd crossings.  Since it is not possible to have self-crossings from more than one component involved in a Reidemeister move we can instead use separate variables $z_i$ for each component $i$ of a link as shown in Figure \ref{fig:GenParityAlexBQ}.  We utilized Mathematica to determine the values for maps for $L_{\{i,j \}}: X\times X \rightarrow X\times X$ which satisfy the definition of a Generalized Link Parity Biquandle pictured diagrammatically as in Figure \ref{fig:GenParityAlexBQ}.\\

Once again we may define a polynomial as above which we refer to as the Generalized Link Parity Alexander Polynomial.  An analogous proof to that of Theorem \ref{thm:main} gives the following:\\

\begin{thm}\label{thm:genodd}
Given a virtual link $K$ with components labeled $1, \ldots, k$. For $i \in \{1, \ldots, k \}$ let $o_{i}$ be the minimum number of odd crossings in component $i$ of any diagram of $K$, and suppose $z_{i}^{e_{max}}$ and $z_{i}^{e_{min}}$ are, respectively, the highest and lowest powers of $z_{i}$ appearing in the Generalized Link Parity Alexander Polynomial of $K$. Then $ \max{(\left| e_{max} \right|,\left| e_{min} \right|)}  \leq o_{i}$.\\
\end{thm}

\begin{thm}\label{thm:genlink}
Given a virtual link $K$ with components labeled $1, \ldots, k$. For $i \neq j \in \{1, \ldots, k \}$ let $l_{\{i,j\}}$ be the minimum number of link crossings between components $i$ and $j$ of any diagram of $K$, and suppose $w_{\{i,j\}}^{e_{max}}$ and $w_{\{i,j\}}^{e_{min}}$ are, respectively, the highest and lowest powers of $w_{\{i,j\}}$ appearing in the Generalized Link Parity Alexander Polynomial of $K$. Then $ \max{(\left| e_{max} \right|,\left| e_{min} \right|)}  \leq l_{\{i,j\}}$.\\
\end{thm}

 \begin{exa} The virtual link in Figure \ref{fig:LinkParityExample} has Generalized Link Parity Alexander Polynomial.\\
 \begin{eqnarray*} & & \frac{s t}{w z_{1}^2}+\frac{s t}{w z_{2}^2}-\frac{s t}{w z_{1}^2 z_{2}^2}-\frac{1}{s t w z_{1}^2 z_{2}^2}-\frac{s t}{w}-\frac{w}{s t}+\frac{1}{s t z_{1}^2}+ \\
 & & \frac{1}{s t z_{2}^2}-\frac{1}{w z_{1}^2}-\frac{1}{w z_{2}^2}+\frac{2}{w z_{1}^2 z_{2}^2}-\frac{1}{z_{1}^2}-\frac{1}{z_{2}^2}+2 \\
\end{eqnarray*}
 Note that Theorems \ref{thm:genodd} and \ref{thm:genlink} prove the diagram in Figure \ref{fig:LinkParityExample} is minimal both in the number of odd crossings in each component as well as in the number of crossing between components.\\
 \end{exa}

\subsection{Extensions to Virtual Crossings}

\begin{figure}[H]
\centering
    \includegraphics[height=2.5cm]{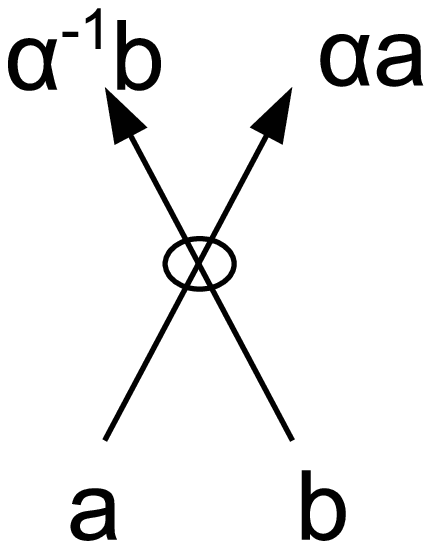}
\caption{}
\label{fig:ManturovTwist}
\end{figure}

\begin{figure}[H]
\centering
    \includegraphics[height=5cm]{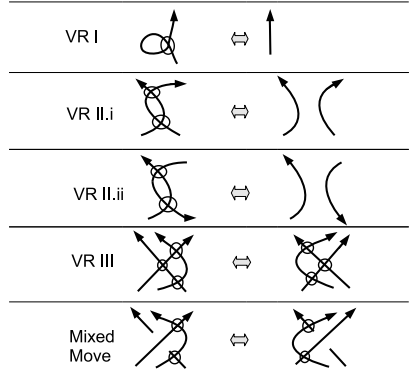}
\caption{Oriented Virtual Reidemeister Moves}
\label{fig:OVRMs}
\end{figure}

Manturov's twist relation for virtual crossings (\cite{Manturov}, \cite{KR1}) for the Generalized Alexander Biquandle as shown in Figure \ref{fig:ManturovTwist} allows us to further extend the Alexander Biquandle and the Generalized Link Parity Alexander Biquandle to what we call the $\alpha$-Alexander Biquandle and $\alpha$-Generalized Link Parity Alexander Biquandle. Denote the relation and associated map from $X \times X \to X \times X$ in Figure \ref{fig:ManturovTwist} by $V$. Notice that replacing $\alpha$ by $w$ in $V$ we have the link parity relation $L$ for the Generalized Link Parity Alexander Biquandle. Thus $V$ satisfies the axioms of $L$. In other words, the $\alpha$-Alexander Biquandle and $\alpha$-Generalized Link Parity Alexander Biquandle are invariant under  both oriented virtual Reidemeister II moves, the oriented virtual Reidemeister III move as well as the oriented Mixed Moves (Figure \ref{fig:OVRMs}). Moreover, it is easy to check that $V$ satisfies Axiom 4 of the Biquandle map $B$ implying the $\alpha$-Alexander Biquandle and $\alpha$-Generalized Link Parity Alexander Biquandle are invariant under the Virtual Reidemeister I Move. Hence the $\alpha$-Alexander Biquandle and $\alpha$-Generalized Link Parity Alexander Biquandle are invariants of virtual knots.\\

Starting with the Generalized Alexander Biquandle incorporating Manturov's Twist creates the $\alpha$-Generalized Alexander Biquandle and the $\alpha$-Generalized Link Parity Alexander Biquandle with respective polynomials we refer to as $\alpha$-Sawollek and $\alpha$-Generalized Link Parity Alexander.  An analogous proof to that of Theorem \ref{thm:main} gives the following:\\

\begin{thm}\label{thm:virtual}
Given a virtual link $K$, let $n_{v}$ be the minimum number of virtual crossings in any diagram of $K$, and suppose $\alpha^{e_{max}}$ and $\alpha^{e_{min}}$ are, respectively, the highest and lowest powers of $\alpha$ appearing in the $\alpha$-Sawollek polynomial (or ($\alpha$-Generalized Link Parity Alexander polynomial) of $K$. Then $ \max{(\left| e_{max} \right|,\left| e_{min} \right|)}  \leq n_{v}$.\\
\end{thm}

 \begin{exa}
 Calculating  the $\alpha$-Sawollek polynomial and the $\alpha$-Generalized Link Parity Alexander polynomial for the virtual knot $3.1$ in Figure \ref{fig:3.1} we get the respectively:
 \[\alpha^{-1}(-s + t^{-1})+(-1 + st) + \alpha (s^{-1} - s^{-2} t^{-1}) +\alpha^{2}(s^{-2} - t s^{-1}) \]
  and
 \[s^{-1} t^{-1} (\alpha^{2} z^{-2} -1)\]
	
It follows from Theorem \ref{thm:main}, Corollary \ref{cor:maincor} and Theorem \ref{thm:virtual} that the diagram in Figure \ref{fig:3.1} is minimal in virtual crossing number, odd crossing number and total crossing number.
\end{exa}

\begin{exa}The virtual link in Figure \ref{fig:LinkParityExample} has $\alpha$-Generalized Link Parity Alexander Polynomial.\\

\begin{eqnarray*}
& & -\frac{s t {\alpha}^5}{z_{2}^2 w z_{1}^2}  -  \frac{{\alpha}^5}{s t z_{2}^2 w z_{1}^2} +  \frac{s t {\alpha}^3} {z_{2}^2 w} +  \frac{{\alpha}^2}{s t z_{2}^2} + \frac{s t {\alpha}^3}{w z_{1}^2} + \frac{{\alpha}^2}{s t z_{1}^2} -  \\
& &  \frac{s t {\alpha}}{w} - \frac{w}{s t {\alpha}} + \frac {2 {\alpha}^5} {z_{2}^2 w z_ {1}^2} - \frac {{\alpha}^3} {z_{2}^2 w} - \frac {{\alpha}^2} {z_{2}^2} - \frac {{\alpha}^3} {w z_{1}^2} - \frac {{\alpha}^2} {z_{1}^2} + 2
\end{eqnarray*}

Thus Theorems \ref{thm:genodd} and \ref{thm:genlink} prove the diagram in Figure \ref{fig:LinkParityExample} is minimal both in the number of odd crossings in each component as well as in the number of crossing between components and Theorem \ref{thm:virtual} shows that the diagram is minimal with respect to virtual crossing number.\\
\end{exa}

\subsection{Further Questions and Remarks}

We have only begun to scratch the surface in our search for parity biquandles. In the linear case, Bartholomew and Fenn have shown in \cite{Quaternionic} there are additional quaternionic biquandles with coefficients in the Hurwitz ring. One would expect to find similar results to the linear biquandle structures studied here.  Additionally, Bartholomew and Fenn  \cite{BFSmall} point out the nonlinear biquandles of Wada \cite{Wada} and Silver and Williams \cite{SilverWilliams}. It is hopeful that additional useful examples will arise from these structures.\\

\appendix
\section{Appendix}

\begin{longtable}{l|p{5.7cm}|p{4.8cm}}
Knot & Sawollek Polynomial & z-Parity Alexander Polynomial \\
\hline \endhead
    2.1 & $\left(s^2-s\right) t^2+\left(1-s^2\right) t+s-1$ & $0$ \\ \hline
    3.1 & $ \frac{1-\frac{1}{s^2}}{t}+\frac{1}{s^2}+\left(s-\frac{1}{s}\right) t-s+\frac{1}{s}-1 $ & $\frac{\frac{1}{s t}-1}{z^2}-\frac{1}{s t}+1$ \\ \hline
    3.2 & $ \frac{\frac{1}{s}-\frac{1}{s^2}}{t^2}+\frac{\frac{1}{s^2}-1}{t}-\frac{1}{s}+1 $ & $ 0$ \\ \hline
    3.3 & $ \frac{\frac{1}{s^3}-1}{t}-\frac{1}{s^2}+\frac{\frac{1}{s^2}-\frac{1}{s^3}}{t^3}+1 $ & $ \frac{\frac{1}{s t}-1}{z^2}-\frac{1}{s t}+1$ \\ \hline
    3.4 & $ \frac{\frac{1}{s}-\frac{1}{s^2}}{t}+\frac{1}{s^2}+\left(1-\frac{1}{s}\right) t-1 $ & $ \frac{1-s t}{z^2}+s t-1$ \\ \hline
    3.5 & $ \frac{\frac{1}{s}-\frac{1}{s^3}}{t^3}+\frac{\frac{1}{s^3}-\frac{1}{s}}{t}+\frac{\frac{1}{s^2}-1}{t^2}-\frac{1}{s^2}+1 $ & $ -\frac{1}{s^3 t^3}+\frac{1}{s^3 t}+\frac{1}{s^2 t^2}-\frac{1}{s^2}+\frac{1}{s t^3}-\frac{1}{s t}-\frac{1}{t^2}+1$ \\ \hline
    3.6 & $ 0 $ & $ 0$ \\ \hline
    3.7 & $ \frac{1-\frac{1}{s^2}}{t^2}+\frac{1}{s^2}+\left(s-\frac{1}{s}\right) t+\frac{\frac{1}{s}-s}{t}-1 $ & $ -\frac{1}{s^2 t^2}+\frac{1}{s^2}+s t-\frac{s}{t}-\frac{t}{s}+\frac{1}{s t}+\frac{1}{t^2}-1$ \\ \hline
 4.1 & $\frac{\frac{2}{s^3}-\frac{2}{s}}{t}-\frac{1}{s^2}+\frac{-\frac{1}{s^4}+\frac{2}{s^3}-\frac{1}{s^2}}{t^4}+\frac{\frac{2}{s^4}-\frac{2}{s^3}-\frac{2}{s^2}+\frac{2}{s}}{t^3}+\frac{-\frac{1}{s^4}-\frac{2}{s^3}+\frac{4}{s^2}-1}{t^2}+1$ & $0$ \\ \hline
4.2 & $\left(-s^2+2 s-\frac{2}{s}+1\right) t+\frac{-\frac{1}{s^2}-2 s+\frac{2}{s}+1}{t}+s^2+\frac{1}{s^2}+(1-s) t^2+\frac{1-\frac{1}{s}}{t^2}+s+\frac{1}{s}-4$ & $0$ \\ \hline
4.3 & $\frac{-\frac{2}{s^3}+\frac{1}{s}+1}{t}+\frac{1}{s^2}+\frac{\frac{1}{s^4}-\frac{1}{s^3}}{t^4}+\frac{-\frac{2}{s^4}+\frac{1}{s^3}+\frac{1}{s^2}}{t^3}+\frac{\frac{1}{s^4}+\frac{2}{s^3}-\frac{2}{s^2}-\frac{1}{s}}{t^2}-1$ & $0$ \\ \hline
4.4 & $\frac{\frac{1}{s^2}-\frac{1}{s}}{t^3}+\frac{2-\frac{2}{s^2}}{t^2}+\frac{\frac{1}{s^2}-s+\frac{2}{s}-2}{t}+s-\frac{1}{s}$ & $0$ \\ \hline
4.5 & $\frac{\frac{1}{s}-\frac{1}{s^3}}{t^2}-\frac{2}{s^2}+\frac{\frac{1}{s^3}+\frac{2}{s^2}-\frac{2}{s}-1}{t}+\left(\frac{1}{s}-1\right) t+2$ & $0$ \\ \hline
4.6 & $\left(s^2-s \right) t-s^2+\frac{\frac{1}{s}-1}{t^2}+\frac{2 s-\frac{1}{s}-1}{t}-s+2$ & $0$ \\ \hline
4.7 & $\frac{\frac{1}{s}-\frac{1}{s^3}}{t^2}+\frac{\frac{1}{s^2}-1}{t}+\frac{\frac{1}{s^3}-\frac{1}{s^4}}{t^4}+\frac{\frac{1}{s^4}-\frac{1}{s^2}}{t^3}-\frac{1}{s}+1$ & $0$ \\ \hline
4.8 & $0$ & $0$ \\ \hline
4.9 & $\frac{\frac{2}{s^3}-\frac{2}{s}}{t^2}+\frac{-\frac{2}{s^2}+\frac{1}{s}+1}{t}+\frac{\frac{1}{s^4}-\frac{1}{s^3}}{t^4}+\frac{-\frac{1}{s^4}-\frac{1}{s^3}+\frac{2}{s^2}}{t^3}+\frac{1}{s}-1$ & $\frac{1}{s^2 t^2}+z \left(\frac{1}{s}-\frac{1}{s^2 t}\right)+\frac{\frac{1}{t}-\frac{1}{s t^2}}{z}-1$ \\ \hline
4.10 & $\frac{\frac{1}{s}-\frac{1}{s^3}}{t^2}-\frac{1}{s^2}+\frac{\frac{1}{s^3}+\frac{1}{s^2}-\frac{1}{s}-1}{t}+1$ & $\frac{\frac{1}{s^2 t}+\frac{1}{s t^2}-\frac{1}{s}-\frac{1}{t}}{z}-\frac{1}{s^2 t^2}+1$ \\ \hline
4.11 & $\frac{\frac{1}{s^2}-\frac{1}{s}}{t^3}-\frac{1}{s^2}+\frac{-\frac{1}{s^3}-\frac{2}{s^2}+\frac{1}{s}+2}{t^2}+\frac{\frac{1}{s^3}+\frac{2}{s^2}-s+\frac{1}{s}-3}{t}+s-\frac{1}{s}+1$ & $\frac{\frac{1}{s}-t}{z}+z \left(t-\frac{1}{s}\right)$ \\ \hline
4.12 & $\left(s^2-\frac{2}{s}+1\right) t+\frac{-\frac{1}{s^2}+2 s-1}{t}-s^2+\frac{1}{s^2}+(1-s) t^2+\frac{\frac{1}{s}-1}{t^2}-s+\frac{1}{s}$ & $0$ \\ \hline
4.13 & $\frac{\frac{1}{s^2}-\frac{1}{s}}{t}-\frac{1}{s^2}+(s-1) t^2+\left(-s+\frac{2}{s}-1\right) t-\frac{1}{s}+2$ & $0$ \\ \hline
4.14 & $\frac{-\frac{1}{s^2}-s+\frac{1}{s}+1}{t}+\frac{1}{s^2}+\frac{1-\frac{1}{s}}{t^2}+\left(s-\frac{1}{s}\right) t+\frac{1}{s}-2$ & $s^2 t^2+z \left(s-s^2 t \right)+\frac{t-s t^2}{z}-1$ \\ \hline
4.15 & $\frac{\frac{1}{s^3}-\frac{1}{s^4}}{t^4}+\frac{\frac{1}{s^3}-\frac{1}{s^2}}{t^3}+\frac{-\frac{1}{s^3}+\frac{2}{s^2}-1}{t}+\frac{\frac{1}{s^4}-\frac{1}{s^3}-\frac{1}{s^2}+\frac{1}{s}}{t^2}-\frac{1}{s}+1$ & $\frac{\frac{1}{s^2 t}+\frac{1}{s t^2}-\frac{1}{s}-\frac{1}{t}}{z}-\frac{1}{s^2 t^2}+1$ \\ \hline
4.16 & $0$ & $\frac{1}{s^2 t^2}+z \left(\frac{1}{s}-\frac{1}{s^2 t}\right)+\frac{\frac{1}{t}-\frac{1}{s t^2}}{z}-1$ \\ \hline
4.17 & $(s-1) t^2+\left(\frac{1}{s}-1\right) t+\frac{1-\frac{1}{s}}{t}-s+1$ & $\frac{\frac{1}{s}-t}{z}+z \left(t-\frac{1}{s}\right)$ \\ \hline
4.18 & $\frac{\frac{1}{s^2}-\frac{1}{s}}{t^2}+\frac{1-\frac{1}{s^2}}{t}+\frac{1}{s}-1$ & $0$ \\ \hline
4.19 & $\frac{1-\frac{1}{s}}{t^2}+(s-1) t+\frac{1-s}{t}+\frac{1}{s}-1$ & $\frac{s-\frac{1}{t}}{z}+z \left(\frac{1}{t}-s \right)$ \\ \hline
4.20 & $\frac{\frac{1}{s^2}-\frac{1}{s}}{t}-\frac{1}{s^2}+\left(\frac{1}{s}-1\right) t+1$ & $\frac{s^2 t+s t^2-s-t}{z}-s^2 t^2+1$ \\ \hline
4.21 & $\left(s^2-s \right) t^3+\left(s^2-2 s+1\right) t^2+\left(1-s^2\right) t-s^2+\frac{s-1}{t}+2 s-1$ & $s^2 t^2+z \left(s-s^2 t \right)+\frac{t-s t^2}{z}-1$ \\ \hline
4.22 & $\left(\frac{1}{s^2}+\frac{1}{s}-2\right) t-\frac{1}{s^2}+\left(s-\frac{1}{s}\right) t^2+\frac{1-\frac{1}{s}}{t}-s+\frac{1}{s}+1$ & $\frac{s^2 t+s t^2-s-t}{z}-s^2 t^2+1$ \\ \hline
4.23 & $\frac{\frac{1}{s^2}-\frac{1}{s}}{t^3}+\frac{\frac{1}{s}-\frac{1}{s^2}}{t}+\frac{1-\frac{1}{s}}{t^2}+\frac{1}{s}-1$ & $z \left(-\frac{1}{s^2 t}-\frac{1}{s t^2}+\frac{1}{s}+\frac{1}{t}\right)+\frac{1}{s^2 t^2}-1$ \\ \hline
4.24 & $\left(s^2-1\right) t^3+\left(s^2-2 s+\frac{1}{s}\right) t^2+\left(-s^2-s+2\right) t-s^2+\frac{s-1}{t}+2 s-\frac{1}{s}$ & $z \left(s^2 (-t)-s t^2+s+t \right)+s^2 t^2-1$ \\ \hline
4.25 & $\frac{\frac{1}{s}-\frac{1}{s^2}}{t^2}+\frac{\frac{1}{s^2}-\frac{1}{s}}{t}+\frac{\frac{1}{s^4}-\frac{1}{s^3}}{t^4}+\frac{\frac{1}{s^3}-\frac{1}{s^4}}{t^3}$ & $0$ \\ \hline
4.26 & $\left(\frac{1}{s}-s^2\right) t^2+\left(s^2-\frac{1}{s^2}+s-1\right) t+\frac{\frac{1}{s^2}-\frac{1}{s}}{t}-s+1$ & $0$ \\ \hline
4.27 & $\frac{\frac{1}{s^3}-\frac{1}{s}}{t^2}+\frac{\frac{1}{s}-\frac{1}{s^2}}{t}+\frac{\frac{1}{s^2}-\frac{1}{s^3}}{t^3}$ & $0$ \\ \hline
4.28 & $\left(1-s^2\right) t^3+\left(s^2-s \right) t+\left(s-\frac{1}{s}\right) t^2+\frac{\frac{1}{s}-1}{t}$ & $0$ \\ \hline
4.29 & $\frac{\frac{2}{s^2}-\frac{2}{s}}{t^3}+\frac{\frac{1}{s^3}-\frac{2}{s^2}+\frac{1}{s}}{t}+\frac{\frac{1}{s^4}-\frac{2}{s^3}+\frac{1}{s^2}}{t^4}+\frac{-\frac{1}{s^4}+\frac{1}{s^3}-\frac{1}{s^2}+1}{t^2}+\frac{1}{s}-1$ & $\frac{-\frac{1}{s^2 t}-\frac{1}{s t^2}+\frac{1}{s}+\frac{1}{t}}{z}+\frac{1}{s^2 t^2}-1$ \\ \hline
4.30 & $\frac{-\frac{1}{s^3}+\frac{2}{s}-1}{t}+\frac{2}{s^2}+\frac{\frac{1}{s^3}-\frac{2}{s^2}+\frac{1}{s}}{t^2}+\left(1-\frac{1}{s}\right) t-\frac{2}{s}$ & $0$ \\ \hline
4.31 & $\frac{\frac{1}{s}-\frac{1}{s^2}}{t^3}+\frac{\frac{1}{s^2}+\frac{1}{s}-2}{t^2}+\frac{s-\frac{2}{s}+1}{t}-s+1$ & $-\frac{1}{s^2 t^2}+z \left(\frac{1}{s^2 t}-\frac{1}{s}\right)+\frac{\frac{1}{s t^2}-\frac{1}{t}}{z}+1$ \\ \hline
4.32 & $\frac{\frac{1}{s^2}-1}{t}-\frac{1}{s^2}+\left(\frac{1}{s}-s \right) t+s-\frac{1}{s}+1$ & $\frac{t-\frac{1}{s}}{z}+z \left(\frac{1}{s}-t \right)$ \\ \hline
4.33 & $\frac{\frac{1}{s}-\frac{1}{s^2}}{t^2}+\frac{\frac{1}{s^2}-1}{t}-\frac{1}{s}+1$ & $0$ \\ \hline
4.34 & $\frac{\frac{1}{s}-\frac{1}{s^2}}{t}+\frac{1}{s^2}+\left(1-\frac{1}{s}\right) t-1$ & $\frac{s^2 (-t)-s t^2+s+t}{z}+s^2 t^2-1$ \\ \hline
4.35 & $\frac{\frac{1}{s}-1}{t^2}+(1-s) t+\frac{s-1}{t}-\frac{1}{s}+1$ & $\frac{\frac{1}{t}-s}{z}+z \left(s-\frac{1}{t}\right)$ \\ \hline
4.36 & $\left(s-s^2\right) t^3+\left(-s^2+2 s-1\right) t^2+\left(s^2-1\right) t+s^2+\frac{1-s}{t}-2 s+1$ & $-s^2 t^2+z \left(s^2 t-s \right)+\frac{s t^2-t}{z}+1$ \\ \hline
4.37 & $\frac{1-\frac{1}{s^4}}{t^2}+\frac{\frac{1}{s^2}-\frac{1}{s}}{t^3}+\frac{\frac{1}{s^4}-\frac{1}{s^3}}{t^4}+\frac{\frac{1}{s^3}-\frac{1}{s^2}}{t}+\frac{1}{s}-1$ & $\frac{1}{s^2 t^2}+\frac{\frac{1}{s}-\frac{1}{s^2 t}}{z}+z \left(\frac{1}{t}-\frac{1}{s t^2}\right)-1$ \\ \hline
4.38 & $\frac{\frac{1}{s}-1}{t^2}+\frac{s-\frac{1}{s}}{t}-s+1$ & $\frac{\frac{1}{t}-s}{z}+z \left(s-\frac{1}{t}\right)$ \\ \hline
4.39 & $\frac{1}{s^3}+\frac{\frac{1}{s}-\frac{1}{s^2}}{t^2}+\left(1-\frac{1}{s^2}\right) t+\frac{1}{s^2}+\frac{-\frac{1}{s^3}+\frac{1}{s^2}+\frac{1}{s}-1}{t}-\frac{2}{s}$ & $\frac{t-\frac{1}{s}}{z}+z \left(\frac{1}{s}-t \right)$ \\ \hline
4.40 & $\frac{\frac{1}{s}-1}{t}+(1-s) t+s-\frac{1}{s}$ & $s^2 t^2+\frac{s-s^2 t}{z}+z \left(t-s t^2\right)-1$ \\ \hline
4.41 & $0$ & $-\frac{1}{s^2 t^2}+z \left(\frac{1}{s^2 t}-\frac{1}{s}\right)+\frac{\frac{1}{s t^2}-\frac{1}{t}}{z}+1$ \\ \hline
4.42 & $\left(\frac{1}{s^2}-\frac{1}{s}\right) t+\frac{\frac{1}{s}-\frac{1}{s^2}}{t}+\left(1-\frac{1}{s}\right) t^2+\frac{1}{s}-1$ & $\frac{t-\frac{1}{s}}{z}+z \left(\frac{1}{s}-t \right)$ \\ \hline
4.43 & $\frac{\frac{2}{s^3}-\frac{2}{s}}{t^2}+\frac{\frac{2}{s}-\frac{2}{s^2}}{t}+\frac{\frac{2}{s^2}-\frac{2}{s^3}}{t^3}$ & $0$ \\ \hline
4.44 & $\frac{\frac{1}{s}-\frac{1}{s^3}}{t^2}+\frac{\frac{1}{s^2}-\frac{1}{s}}{t}+\frac{\frac{1}{s^3}-\frac{1}{s^2}}{t^3}$ & $0$ \\ \hline
4.45 & $-\frac{1}{s^3}+\frac{\frac{1}{s}-\frac{1}{s^2}}{t^2}+\left(\frac{1}{s^2}-s \right) t+\frac{\frac{1}{s^2}-1}{t}+\frac{\frac{1}{s^3}-\frac{1}{s^2}}{t^3}+s-\frac{1}{s}+1$ & $0$ \\ \hline
4.46 & $\frac{\frac{1}{s^2}-\frac{1}{s}}{t^2}+\frac{-\frac{1}{s^2}-\frac{1}{s}+2}{t}+(s-1) t-s+\frac{2}{s}-1$ & $0$ \\ \hline
4.47 & $\left(1-\frac{1}{s^2}\right) t+\frac{\frac{1}{s^2}-1}{t}+\left(\frac{1}{s}-s \right) t^2+s-\frac{1}{s}$ & $0$ \\ \hline
4.48 & $\frac{\frac{1}{s^3}-\frac{1}{s}}{t}-\frac{1}{s^2}+\frac{\frac{1}{s^3}-\frac{1}{s^4}}{t^4}+\frac{-\frac{1}{s^3}+\frac{2}{s^2}-1}{t^2}+\frac{\frac{1}{s^4}-\frac{1}{s^3}-\frac{1}{s^2}+\frac{1}{s}}{t^3}+1$ & $-\frac{1}{s^2 t^2}+\frac{\frac{1}{s^2 t}-\frac{1}{s}}{z}+z \left(\frac{1}{s t^2}-\frac{1}{t}\right)+1$ \\ \hline
4.49 & $\frac{\frac{1}{s}-\frac{1}{s^2}}{t^3}+\frac{\frac{1}{s^2}-1}{t^2}+\frac{1-\frac{1}{s}}{t}$ & $\frac{s-\frac{1}{t}}{z}+z \left(\frac{1}{t}-s \right)$ \\ \hline
4.50 & $\frac{\frac{1}{s^3}-\frac{1}{s}}{t^2}+\frac{1}{s^2}+\frac{-\frac{1}{s^3}-\frac{1}{s^2}+\frac{1}{s}+1}{t}-1$ & $\frac{-\frac{1}{s^2 t}-\frac{1}{s t^2}+\frac{1}{s}+\frac{1}{t}}{z}+\frac{1}{s^2 t^2}-1$ \\ \hline
4.51 & $\left(s-s^2\right) t+s^2+\frac{1-\frac{1}{s}}{t^2}+\frac{-2 s+\frac{1}{s}+1}{t}+s-2$ & $0$ \\ \hline
4.52 & $\frac{1-\frac{1}{s}}{t}+(s-1) t-s+\frac{1}{s}$ & $-s^2 t^2+\frac{s^2 t-s}{z}+z \left(s t^2-t \right)+1$ \\ \hline
4.53 & $\frac{\frac{1}{s^3}-\frac{1}{s}}{t^2}+\frac{1-\frac{1}{s^2}}{t}+\frac{\frac{1}{s^4}-\frac{1}{s^3}}{t^4}+\frac{\frac{1}{s^2}-\frac{1}{s^4}}{t^3}+\frac{1}{s}-1$ & $0$ \\ \hline
4.54 & $\frac{\frac{1}{s^2}-\frac{1}{s}}{t^3}+\frac{2-\frac{2}{s^2}}{t^2}+\frac{\frac{1}{s^2}-s+\frac{2}{s}-2}{t}+s-\frac{1}{s}$ & $0$ \\ \hline
4.55 & $0$ & $0$ \\ \hline
4.56 & $0$ & $0$ \\ \hline
4.57 & $\frac{\frac{1}{s}-\frac{1}{s^2}}{t^2}+\frac{1}{s^2}+\frac{\frac{1}{s}-1}{t}+\left(1-\frac{1}{s}\right) t-\frac{1}{s}$ & $\frac{\frac{1}{s^2 t}+\frac{1}{s t^2}-\frac{1}{s}-\frac{1}{t}}{z}-\frac{1}{s^2 t^2}+1$ \\ \hline
4.58 & $0$ & $0$ \\ \hline
4.59 & $0$ & $0$ \\ \hline
4.60 & $\frac{\frac{1}{s}-1}{t}+(1-s) t+s-\frac{1}{s}$ & $s^2 t^2+\frac{s-s^2 t}{z}+z \left(t-s t^2\right)-1$ \\ \hline
4.61 & $\frac{\frac{2}{s^3}-\frac{2}{s}}{t^2}+\frac{-\frac{2}{s^2}+\frac{1}{s}+1}{t}+\frac{\frac{1}{s^4}-\frac{1}{s^3}}{t^4}+\frac{-\frac{1}{s^4}-\frac{1}{s^3}+\frac{2}{s^2}}{t^3}+\frac{1}{s}-1$ & $\frac{1}{s^2 t^2}+z \left(\frac{1}{s}-\frac{1}{s^2 t}\right)+\frac{\frac{1}{t}-\frac{1}{s t^2}}{z}-1$ \\ \hline
4.62 & $\frac{1}{s^3}+\left(1-\frac{1}{s^2}\right) t+\frac{\frac{3}{s^2}-s+\frac{1}{s}-3}{t}+\frac{-\frac{1}{s^3}-\frac{2}{s^2}+\frac{2}{s}+1}{t^2}+s-\frac{3}{s}+1$ & $\frac{s-\frac{1}{t}}{z}+z \left(\frac{1}{t}-s \right)$ \\ \hline
4.63 & $-\frac{1}{s^2}+\frac{-\frac{1}{s^3}-\frac{1}{s^2}+\frac{2}{s}}{t^2}+\frac{\frac{1}{s^3}+\frac{2}{s^2}-\frac{1}{s}-2}{t}-\frac{1}{s}+2$ & $\frac{\frac{1}{s}-t}{z}+z \left(t-\frac{1}{s}\right)$ \\ \hline
4.64 & $\left(s^2-\frac{1}{s}\right) t+\frac{s-\frac{1}{s^2}}{t}-s^2+\frac{1}{s^2}-s+\frac{1}{s}$ & $s^2 t^2+z \left(s-s^2 t \right)+\frac{t-s t^2}{z}-1$ \\ \hline
4.65 & $\frac{\frac{1}{s}-\frac{1}{s^2}}{t^3}+\frac{-\frac{1}{s^2}+\frac{2}{s}-1}{t^2}+\frac{\frac{1}{s^2}-1}{t}+\frac{1}{s^2}+\left(1-\frac{1}{s}\right) t-\frac{2}{s}+1$ & $-\frac{1}{s^2 t^2}+\frac{\frac{1}{s^2 t}-\frac{1}{s}}{z}+z \left(\frac{1}{s t^2}-\frac{1}{t}\right)+1$ \\ \hline
4.66 & $\left(\frac{1}{s^2}+s-2\right) t+\frac{-\frac{1}{s^2}-s+2}{t}+\left(s-\frac{1}{s}\right) t^2+\frac{1-\frac{1}{s}}{t^2}-s+\frac{2}{s}-1$ & $\frac{s-\frac{1}{t}}{z}+z \left(\frac{1}{t}-s \right)$ \\ \hline
4.67 & $(s-1) t^2+\left(\frac{1}{s}-1\right) t+\frac{1-\frac{1}{s}}{t}-s+1$ & $\frac{\frac{1}{s}-t}{z}+z \left(t-\frac{1}{s}\right)$ \\ \hline
4.68 & $0$ & $-s^2 t^2+\frac{s^2 t-s}{z}+z \left(s t^2-t \right)+1$ \\ \hline
4.69 & $\frac{\frac{2}{s}-\frac{2}{s^3}}{t^2}+\frac{\frac{2}{s^2}-\frac{1}{s}-1}{t}+\frac{\frac{1}{s^3}-\frac{1}{s^4}}{t^4}+\frac{\frac{1}{s^4}+\frac{1}{s^3}-\frac{2}{s^2}}{t^3}-\frac{1}{s}+1$ & $-\frac{1}{s^2 t^2}+z \left(\frac{1}{s^2 t}-\frac{1}{s}\right)+\frac{\frac{1}{s t^2}-\frac{1}{t}}{z}+1$ \\ \hline
4.70 & $\frac{\frac{1}{s^3}-\frac{1}{s}}{t^2}+\frac{1}{s^2}+\frac{-\frac{1}{s^3}-\frac{1}{s^2}+\frac{1}{s}+1}{t}-1$ & $\frac{-\frac{1}{s^2 t}-\frac{1}{s t^2}+\frac{1}{s}+\frac{1}{t}}{z}+\frac{1}{s^2 t^2}-1$ \\ \hline
4.71 & $0$ & $0$ \\ \hline
4.72 & $0$ & $0$ \\ \hline
4.73 & $\frac{1-\frac{1}{s^4}}{t^2}+\frac{\frac{2}{s^2}-\frac{2}{s}}{t^3}-\frac{1}{s^2}+\frac{\frac{2}{s^3}-\frac{2}{s^2}}{t}+\frac{\frac{1}{s^4}-\frac{2}{s^3}+\frac{1}{s^2}}{t^4}+\frac{2}{s}-1$ & $0$ \\ \hline
4.74 & $\frac{\frac{1}{s}-\frac{1}{s^2}}{t^3}+\frac{\frac{1}{s^2}+s-\frac{2}{s}}{t}+\frac{\frac{2}{s}-2}{t^2}-s-\frac{1}{s}+2$ & $0$ \\ \hline
4.75 & $0$ & $0$ \\ \hline
4.76 & $0$ & $0$ \\ \hline
4.77 & $0$ & $0$ \\ \hline
4.78 & $\frac{\frac{1}{s^4}-\frac{1}{s}}{t}-\frac{1}{s^3}+\frac{\frac{1}{s}-\frac{1}{s^2}}{t^3}+\frac{\frac{1}{s^2}-1}{t^2}+\frac{\frac{1}{s^3}-\frac{1}{s^4}}{t^4}+1$ & $z \left(\frac{1}{s^2 t}+\frac{1}{s t^2}-\frac{1}{s}-\frac{1}{t}\right)-\frac{1}{s^2 t^2}+1$ \\ \hline
4.79 & $\frac{1}{s^3}+\left(\frac{1}{s}-\frac{1}{s^2}\right) t+\frac{\frac{1}{s^2}-\frac{1}{s}}{t}-\frac{1}{s^2}+\frac{\frac{1}{s^2}-\frac{1}{s^3}}{t^2}$ & $\frac{-\frac{1}{s^2 t}-\frac{1}{s t^2}+\frac{1}{s}+\frac{1}{t}}{z}+\frac{1}{s^2 t^2}-1$ \\ \hline
4.80 & $\frac{1-\frac{1}{s^4}}{t}+\frac{1}{s^3}+\frac{\frac{1}{s^4}-\frac{1}{s^3}}{t^4}-1$ & $0$ \\ \hline
4.81 & $-\frac{1}{s^3}+\left(\frac{1}{s^2}-1\right) t+\frac{\frac{1}{s^3}-\frac{1}{s^2}}{t^2}+1$ & $0$ \\ \hline
4.82 & $\frac{\frac{1}{s}-\frac{1}{s^2}}{t^3}+\frac{\frac{1}{s^4}-\frac{1}{s^2}}{t^4}+\frac{\frac{1}{s^3}+\frac{1}{s^2}-\frac{2}{s}}{t}+\frac{-\frac{1}{s^4}-\frac{1}{s^3}+\frac{1}{s^2}+\frac{1}{s}}{t^2}$ & $z^2 \left(\frac{1}{s t}-\frac{1}{s^2 t^2}\right)+\frac{\frac{1}{s t}-\frac{1}{s^2 t^2}}{z^2}+\frac{2}{s^2 t^2}-\frac{2}{s t}$ \\ \hline
4.83 & $\frac{1}{s^3}+\frac{1-\frac{1}{s^2}}{t^2}+\left(s-\frac{1}{s^2}\right) t+\frac{1}{s^2}+\frac{-\frac{1}{s^3}+\frac{1}{s^2}-s+\frac{1}{s}}{t}-\frac{1}{s}-1$ & $z^2 (1-s t)+\frac{1-\frac{1}{s t}}{z^2}+s t+\frac{1}{s t}-2$ \\ \hline
4.84 & $\frac{\frac{1}{s^2}-1}{t^2}-\frac{1}{s^2}+\left(\frac{1}{s}-1\right) t+\frac{s-\frac{2}{s}+1}{t}-s+\frac{1}{s}+1$ & $z^2 \left(\frac{1}{s t}-\frac{1}{s^2 t^2}\right)+\frac{\frac{1}{s t}-\frac{1}{s^2 t^2}}{z^2}+\frac{2}{s^2 t^2}-\frac{2}{s t}$ \\ \hline
4.85 & $\frac{\frac{1}{s^2}-1}{t^2}-\frac{1}{s^2}+\left(\frac{1}{s}-s \right) t+\frac{s-\frac{1}{s}}{t}+1$ & $\frac{1}{s^2 t^2}-\frac{1}{s^2}+s (-t)+\frac{s}{t}+\frac{t}{s}-\frac{1}{s t}-\frac{1}{t^2}+1$ \\ \hline
4.86 & $\frac{1-\frac{1}{s^2}}{t^2}+\frac{1}{s^2}+\left(s-\frac{1}{s}\right) t+\frac{\frac{1}{s}-s}{t}-1$ & $-\frac{1}{s^2 t^2}+\frac{1}{s^2}+s t-\frac{s}{t}-\frac{t}{s}+\frac{1}{s t}+\frac{1}{t^2}-1$ \\ \hline
4.87 & $\frac{\frac{1}{s^3}-1}{t^2}-\frac{1}{s^3}+\frac{\frac{1}{s^2}-\frac{1}{s^4}}{t^4}+\frac{\frac{1}{s^4}-\frac{1}{s^2}}{t}+1$ & $z^2 \left(\frac{1}{s^2 t^2}-1\right)-\frac{1}{s^2 t^2}+1$ \\ \hline
4.88 & $\left(1-\frac{1}{s^2}\right) t+\frac{1}{s^2}+\left(\frac{1}{s}-1\right) t^2-\frac{1}{s}$ & $\frac{\frac{1}{s^2 t^2}-1}{z^2}-\frac{1}{s^2 t^2}+1$ \\ \hline
4.89 & $\frac{1-\frac{1}{s^4}}{t^2}+\frac{1}{s^2}+\frac{\frac{1}{s^4}-\frac{1}{s^2}}{t^4}-1$ & $\frac{1}{s^4 t^4}-\frac{1}{s^4 t^2}-\frac{1}{s^2 t^4}+\frac{1}{s^2}+\frac{1}{t^2}-1$ \\ \hline
4.90 & $0$ & $0$ \\ \hline
4.91 & $\frac{\frac{1}{s}-\frac{1}{s^4}}{t^4}+\frac{\frac{1}{s^3}-\frac{1}{s}}{t^2}-\frac{1}{s^3}+\frac{\frac{1}{s^2}-1}{t^3}+\frac{\frac{1}{s^4}-\frac{1}{s^2}}{t}+1$ & $0$ \\ \hline
4.92 & $\frac{\frac{1}{s^4}-\frac{1}{s}}{t^4}+\frac{\frac{1}{s}-\frac{1}{s^4}}{t}+\frac{1-\frac{1}{s^3}}{t^3}+\frac{1}{s^3}-1$ & $0$ \\ \hline
4.93 & $\frac{\frac{1}{s^3}-1}{t}-\frac{1}{s^3}+\frac{\frac{1}{s}-\frac{1}{s^2}}{t^2}+\left(\frac{1}{s^2}-\frac{1}{s}\right) t+1$ & $\frac{\frac{1}{s^2 t^2}-1}{z^2}-\frac{1}{s^2 t^2}+1$ \\ \hline
4.94 & $\frac{\frac{1}{s^3}-\frac{1}{s}}{t^2}+\frac{\frac{1}{s}-\frac{1}{s^2}}{t}+\frac{\frac{1}{s^2}-\frac{1}{s^3}}{t^3}$ & $0$ \\ \hline
4.95 & $\frac{\frac{1}{s^3}-1}{t^3}-\frac{1}{s^3}+\frac{s-\frac{1}{s^2}}{t^2}+\left(\frac{1}{s^2}-s \right) t+1$ & $0$ \\ \hline
4.96 & $\frac{\frac{1}{s^3}-\frac{1}{s}}{t^3}+\frac{\frac{1}{s}-\frac{1}{s^3}}{t}+\frac{1-\frac{1}{s^2}}{t^2}+\frac{1}{s^2}-1$ & $z^2 \left(\frac{1}{s t}-\frac{1}{s^2 t^2}\right)+\frac{\frac{1}{s t}-\frac{1}{s^2 t^2}}{z^2}+\frac{2}{s^2 t^2}-\frac{2}{s t}$ \\ \hline
4.97 & $\frac{\frac{1}{s^2}-s}{t^2}+\frac{s^2+s-\frac{2}{s}}{t}-s^2-\frac{1}{s^2}+\left(\frac{1}{s}-1\right) t+\frac{1}{s}+1$ & $z^2 (1-s t)+\frac{1-\frac{1}{s t}}{z^2}+s t+\frac{1}{s t}-2$ \\ \hline
4.98 & $0$ & $0$ \\ \hline
4.99 & $0$ & $0$ \\ \hline
4.100 & $\frac{\frac{1}{s^3}-\frac{1}{s}}{t^2}+\frac{1-\frac{1}{s^2}}{t}+\frac{\frac{1}{s^4}-\frac{1}{s^3}}{t^4}+\frac{\frac{1}{s^2}-\frac{1}{s^4}}{t^3}+\frac{1}{s}-1$ & $0$ \\ \hline
4.101 & $\frac{1}{s^3}+\frac{1-\frac{1}{s^2}}{t^3}+\left(1-\frac{1}{s^2}\right) t+\frac{\frac{3}{s^2}-3}{t}+\frac{-\frac{1}{s^3}-\frac{1}{s^2}-s+\frac{3}{s}}{t^2}+s-\frac{3}{s}+1$ & $0$ \\ \hline
4.102 & $\left(s^2+\frac{1}{s^2}-2\right) t+\frac{-s^2-\frac{1}{s^2}+2}{t}+\left(s-\frac{1}{s}\right) t^2+\frac{s-\frac{1}{s}}{t^2}-2 s+\frac{2}{s}$ & $0$ \\ \hline
4.103 & $\frac{\frac{1}{s^3}-1}{t^3}+\frac{s-\frac{1}{s^2}}{t^2}+\frac{\frac{1}{s^2}-\frac{1}{s^3}}{t}+\left(\frac{1}{s}-s \right) t-\frac{1}{s}+1$ & $\frac{1-\frac{1}{s^2 t^2}}{z^2}+\frac{1}{s^2 t^2}-1$ \\ \hline
4.104 & $\left(\frac{1}{s}-s^2\right) t^2+\frac{\frac{1}{s^2}-s}{t^2}+\left(s-\frac{1}{s^2}\right) t+\frac{s^2-\frac{1}{s}}{t}$ & $0$ \\ \hline
4.105 & $0$ & $0$ \\ \hline
4.106 & $\frac{\frac{1}{s}-\frac{1}{s^3}}{t^3}+\frac{\frac{1}{s^3}-\frac{1}{s}}{t}+\frac{\frac{1}{s^2}-1}{t^2}-\frac{1}{s^2}+1$ & $-\frac{1}{s^3 t^3}+\frac{1}{s^3 t}+\frac{1}{s^2 t^2}-\frac{1}{s^2}+\frac{1}{s t^3}-\frac{1}{s t}-\frac{1}{t^2}+1$ \\ \hline
4.107 & $\left(1-s^2\right) t^2+\frac{1-\frac{1}{s^2}}{t^2}+s^2+\frac{1}{s^2}+\left(2 s-\frac{2}{s}\right) t+\frac{\frac{2}{s}-2 s}{t}-2$ & $s^2 \left(-t^2\right)-\frac{1}{s^2 t^2}+s^2+\frac{1}{s^2}+2 s t-\frac{2 s}{t}-\frac{2 t}{s}+\frac{2}{s t}+t^2+\frac{1}{t^2}-2$ \\ \hline
4.108 & $0$ & $0$ \\ \hline
  \end{longtable}

\end{document}